\documentclass[twoside,a4paper,12pt,centertags]{amsart}
\usepackage{amsmath,amssymb,verbatim,vmargin}
\usepackage[bookmarks=true]{hyperref}   % Hyperref in DVI and PDF

\theoremstyle{plain}
\newtheorem{thm}{Theorem}[section]
\newtheorem{lem}[thm]{Lemma}
\newtheorem{cor}[thm]{Corollary}
\newtheorem{prop}[thm]{Proposition}

\theoremstyle{definition}
\newtheorem{defn}[thm]{Definition}
\newtheorem{rem}[thm]{Remark}

\title{The Kato square root problem for mixed boundary
value problems}
\newcommand{\Mcc}{{M\raise.55ex\hbox{\lowercase{c}}}}
\author{Andreas Axelsson} \author{Stephen Keith} \author{Alan \Mcc Intosh}
\address{Centre for Mathematics and its Applications\\
Australian National University \\ Canberra, ACT 0200, Australia}
\email{andreas.axelsson@math.u-psud.fr}
\email{stephen.keith@anu.edu.au}
\email{Alan.McIntosh@maths.anu.edu.au}
\subjclass{35J25, 47B44, 47F05}
\newcommand{\qedend}{}

\newcommand{\semic}{{:}}       %FIX

\newcommand{\rnum}{{\mathbf R}}
\newcommand{\cnum}{{\mathbf C}}
\newcommand{\nnum}{{\mathbf N}}

\newcommand{\mH}{{\mathcal H}}

\newcommand{\mE}{{\mathcal E}}

\newcommand{\mL}{{\mathcal L}}

\newcommand{\mV}{{\mathcal V}}
\newcommand{\mF}{{\mathcal F}}

\DeclareMathOperator{\re}{Re}

\newcommand{\sett}[2]{ \{ #1 \, \semic \, #2 \} }

\newcommand{\brac}[1]{\langle #1 \rangle}
\newcommand{\supp}{\text{{\rm supp}}\,}
\newcommand{\dist}{\text{{\rm dist}}\,}

\newcommand{\nul}{\textsf{N}}
\newcommand{\hol}{\textsf{H}}
\newcommand{\ran}{\textsf{R}}
\newcommand{\dom}{\textsf{D}}

\newcommand{\clos}[1]{\overline{#1}}

\newcommand{\dyadic}{\triangle}

\newcommand{\barint}{\mbox{$ave \int$}}

\newcommand{\divv}{{\text{div}}\,}

\def\barint_#1{\mathchoice
            {\mathop{\vrule width 6pt
height 3 pt depth -2.5pt
                    \kern -8.8pt
\intop}\nolimits_{#1}}%
            {\mathop{\vrule width 5pt height
3 pt depth -2.6pt
                    \kern -6.5pt
\intop}\nolimits_{#1}}%
            {\mathop{\vrule width 5pt height
3 pt depth -2.6pt
                    \kern -6pt
\intop}\nolimits_{#1}}%
            {\mathop{\vrule width 5pt height
3 pt depth -2.6pt
          \kern -6pt \intop}\nolimits_{#1}}}

\newcommand{\C}{\cnum}
\newcommand{\R}{\rnum}
\newcommand{\N}{\nnum}

\DeclareMathOperator{\I}{I}
\newcommand{\mC}{{\mathcal C}}

\begin{document}

%\tableofcontents
\maketitle

%
%  SECTION 1 ##############################################################################
%
%
\section{Introduction}

The Kato square root problem for elliptic operators on Lipschitz
domains with mixed boundary conditions can be formulated as
follows. Let $\Omega \subset \R^n$, $n \in \N$, be a Lipschitz
domain, let $\Sigma_1$ be an open subset of the boundary $\Sigma$ of $\Omega$, and define
\begin{equation}  \label{vspace}
 \mV = \left\{ u \in H^1 (\Omega ; \C) :
\supp(\gamma u) \subset \clos{\Sigma_1}  \right\}
\end{equation} 
where $\gamma $ is the trace operator from the Sobolev space
$H^1 (\Omega ; \C)$ to the boundary Sobolev space $H ^{1/2} (\Sigma ; \C)$. 

Given a matrix valued function $A=(a_{jk})$ where $a_{jk} \in
L_\infty(\Omega;\C)$ for each $j,k=0,1,\dots,n$, let $J_A: \mV
\times \mV \longrightarrow \C$ be given by
\begin{equation}  \begin{split}
J_A[u,v] &= \int_\Omega\sum_{j,k=1}^n \left( 
a_{jk} \frac{\partial u}{\partial x_k}\frac{\partial \overline v}{ \partial x_j} + a_{j0} u \frac{
\partial \overline v}{\partial x_j} +
a_{0k}
\frac{\partial u}{ \partial x_k} \overline v +
a_{00} u \overline v \right) \, dx
\end{split}
\end{equation} for every $u,v \in \mV$. 

Suppose that $J_A$ satisfies the following coercivity condition:
there exists $\kappa>0$ such that
\begin{equation} \label{JAU}
\re J_A[u,u] \ge \kappa \left ( \|\nabla u \|^2 + \|u \|^2
\right)\end{equation} for every $ u \in \mV$. Here and below
$(\cdot, \cdot)$ and $\| \cdot \|$ denote the inner product and
norm on $L_2(\Omega;\C)$. Then $J_A$ is a densely defined, closed,
accretive sesquilinear form. Consequently, there exists an
operator $L_A$ on $L_2(\Omega;\C)$ with $\dom(L_A) \subset \mV$
uniquely determined by the property that it is maximal accretive and
satisfies $ J_A[u,v] = (L_A u,v)$ for every $u \in \dom(L_A)$ and
$v \in \mV$. Indeed, $L_A$ is the divergence form operator
$L_Au= -\sum\frac{\partial}{\partial x_j} (a_{jk}\frac{\partial u}{\partial x_k}) -
\sum \frac {\partial}{\partial x_j}(a_{j0}u)+\sum a_{0k}\frac{\partial u}{\partial x_k}
+a_{00}u $ with Dirichlet boundary condition $u=0$ on $\Sigma\setminus
\clos{\Sigma_1}$ and 
natural boundary condition $\sum \nu_j a_{jk} \frac{\partial u}{\partial x_k}
+\sum \nu_j a_{j0} u=0$ on $\Sigma_1$,
defined in an appropriate weak sense.

The square root $\sqrt{L_A}$ of $L_A$ is the
unique maximal accretive operator with 
$\left(\sqrt{L_A}\right)^2=L_A$. For an explanation of the terminology and 
results see \cite[VI -- Theorem 2.1, V -- Theorem 3.35, VI
-- Remark 2.29]{kato}. Also see \cite[Chapter II]{lions}, \cite[Chapter 1]{necas}
and \cite{pryde} for specific material
on forms such as $J_A$ and further references to mixed boundary value problems.
The Kato square root problem is to
determine whether the domain $\dom \left(\sqrt {L_A} \right) =\mV$.

The Kato square root problem for second order elliptic
operators on $\Omega=\R^n$ was solved in \cite{AHLMcT} by 
P.\ Auscher, S.\ Hofmann, M.\ Lacey, A. \Mcc Intosh and Ph.\ Tchamitchian, 
and for higher order elliptic operators and systems
on $\R^n$ in \cite{AHMcT} by Auscher, Hofmann, \Mcc Intosh and Tchamitchian.
The Kato square root problem for second order elliptic operators
on strongly Lipschitz domains with Dirichlet or Neumann boundary
conditions was solved by Auscher and Tchamitchian \cite{Ausc:mix}
who reduced the problem to the Kato problem on $\R^n$ by using
extension maps. As different extensions are required for the 
Dirichlet and Neumann problems, their procedure does not work
for mixed boundary value problems.

The following theorem solves the Kato square root problem for
second order elliptic operators on Lipschitz domains with mixed
boundary conditions. This result is new for both smooth and
Lipschitz domains, and answers a question posed by J.-L. Lions in
1962 \cite[Remark 6.1]{lion:conj}.
We remark that in the case when the coefficients are 
H\"older continuous, the Kato square root problem with mixed
boundary conditions was solved in \cite{Mc2}.

\begin{thm} \label{lions}
Let $\Omega' \subset \R^n$, $n\in\N$, be a smooth domain
which coincides with either the empty set $\emptyset$, 
the half space $\R^n_+$ or $\R^n$ on the complement of a bounded set.
Let $\Sigma_1'\subset\Sigma'=\partial\Omega'$, be a smooth open set, 
which coincides with either the empty set $\emptyset$, the half space 
$\R^{n-1}_+\subset\R^{n-1}$ or $\R^{n-1}$ on the complement 
of a bounded set.

Let $\Omega \subset \R^n$ be a bi-lipschitz image of $\Omega'$
and let $\Sigma_1\subset\Sigma=\partial\Omega$ be the 
corresponding bi-lipschitz image of $\Sigma_1'$.

Define $\mV$, $A$, $J_A$ with the properties specified in 
(\ref{vspace}-\ref{JAU}), and let $L_A$ be the associated maximal accretive operator.

Then $\dom \left({\sqrt{L_A}} \right)= \mV$ with $ \| \sqrt{L_A}
u \| \approx \| \nabla u \| + \|u \|$ for every $u \in \mV.$ The
comparability constant implicit in the use of ``$\approx$''
depends on $\|A\|_{\infty}$ and $\kappa$, 
as well as the constants implicit in the assumptions on $\Omega$
and  $\Sigma_1$. 
\end{thm}

Indeed, a somewhat more general version is presented in Section \ref{slions},
Theorem \ref{local}, concerning elliptic systems with local boundary conditions.
This constitutes an application of results (Theorems \ref{mainthm} and 
Corollary \ref{cortomainthm22}) on homogenous first
order systems $\Gamma$ acting on $L_2(\Omega,\cnum^N)$ which
satisfy $\Gamma^2=0$. We let $\Pi=\Gamma+\Gamma^*$, and consider
perturbations of the type $\Pi_B=\Gamma+B_1\Gamma^*B_2$ where
$B_1$ has positive real part on the range of $\Gamma^*$, $B_2$ has
positive real part on the range of $\Gamma$, and
$\Gamma^*B_2B_1\Gamma^*=0$  and $\Gamma B_1B_2\Gamma =0$. It is
 shown under certain hypotheses that $\Pi_B$ satisfies
quadratic estimates in $L_2(\Omega;\cnum^N)$, and hence that the estimate
$\|\sqrt{{\Pi_B}^2}u\| \approx \|\Pi_Bu\|$ holds.

Techniques developed in the current paper build upon ideas
introduced  by the authors in \cite{AKMC}, where we prove
quadratic estimates for complex perturbations of Dirac-type
operators on $\R^n$ and show that such operators have a bounded
functional calculus. This paper was in turn inspired by the proof of the Kato square root
in \cite{AHLMcT}.
The key idea employed from \cite{AKMC} 
is our utilization of only the first order structure of the 
operator, and subsequent exploitation of the algebra involved 
in the  Hodge decomposition of the first order 
system. Duplicated arguments from \cite{AKMC} have been omitted, 
so the reader is advised to keep a copy of that paper handy.

\subsection{Acknowledgments}

This research was mostly undertaken at the Centre for Mathematics
and its Applications at the Australian National University, and
was supported by the Australian Research Council. The second
author held a visiting position at the School of Mathematics at
the University of New South Wales during the final preparation of
this paper, and thanks them for their hospitality.

%
%  SECTION 2 ##############################################################################
%
%
\section{Quadratic estimates for perturbed Dirac operators}  \label{statement}

In this section we expand on the comments made in the introduction concerning first order elliptic systems.

For an unbounded linear operator $T: \dom(T)\longrightarrow \mH_2$
from a domain $\dom(T)$ in a Hilbert space $\mH_1$ to another
Hilbert space $\mH_2$, denote its null space by $\nul(T)$ and
its range by $\ran(T)$. The operator $T$ is said to be closed when
its graph is a closed subspace of $\mH_1\times \mH_2$. The space
of all bounded linear operators from $\mH_1$ to $\mH_2$ is denoted
$\mL(\mH_1,\mH_2)$, while  $\mL(\mH)=\mL(\mH,\mH)$. See for
example \cite{kato} for more details.

Consider three operators $\{\Gamma, B_1,B_2 \}$ in a Hilbert space
$\mH$ with the following properties.
\begin{itemize}
\item[(H1)] The operator $\Gamma:\dom(\Gamma)\longrightarrow\mH$
is a {\em nilpotent} operator from $\dom(\Gamma)\subset \mH$ to
$\mH$, by which we mean  $\Gamma$ is closed, densely defined and
$\ran(\Gamma)\subset\nul(\Gamma)$. In particular, $\Gamma^2=0$ on
$\dom(\Gamma)$. 
\item[(H2)] The operators $B_1$, $B_2:
\mH\longrightarrow\mH$ are bounded linear operators satisfying the following
accretivity conditions for some $\kappa_1, \kappa_2>0$:
$$   \re (B_1 u,u) \geq \kappa_1 \|u\|^2\quad \text{for all} \quad u
\in\ran(\Gamma^*),
$$
$$  \re (B_2 u,u) \geq \kappa_2 \|u\|^2\quad \text{for all} \quad
u\in\ran(\Gamma).
$$
Let the angles of accretivity be
$$   \omega_1 = \sup_{u\in\ran(\Gamma^*)\setminus\{0\}}|\arg(B_1
u,u)| <\tfrac\pi2\, ,$$
$$
     \omega_2 = \sup_{u\in\ran(\Gamma)\setminus\{0\}}|\arg(B_2 u
,u)| <\tfrac\pi2\,,$$ and set $ \omega
=\tfrac12(\omega_1+\omega_2).$ \item[(H3)] The operators satisfy
$\Gamma^*B_2B_1\Gamma^*=0$ on $\dom(\Gamma^*)$ and $\Gamma
B_1B_2\Gamma =0$ on $\dom(\Gamma)$; that is,  $B_2B_1:
\ran(\Gamma^*) \longrightarrow \nul(\Gamma^*)$ and $B_1B_2:
\ran(\Gamma) \longrightarrow \nul(\Gamma)$. 
\end{itemize}

\begin{defn} Let $\Pi = \Gamma + \Gamma^*$. Also let $\Gamma^*_B = B_1
\Gamma^* B_2 $ and let $\Pi_B= \Gamma + \Gamma^* _B $.
\end{defn}

\begin{defn} \label{psidef} Given $ 0 \leq\mu< 
\frac\pi2$, define the
     closed sectors and
double sector in the complex plane by
\begin{align*}
    S_{\mu+} &= \sett{z\in\cnum}{|\arg z|\le\mu}\cup\{0\}\,,
    \quad S_{\mu-}=-S_{\mu+}\,,\\
    S_{\mu} &= S_{\mu+}\cup S_{\mu-}\,.
\end{align*}\end{defn}

We now summarize consequences of the above hypotheses, proved in
Section 4 of \cite{AKMC}. 
The operator $\Gamma^*_B$ is nilpotent, the operator $\Pi_B$ is closed and densely defined, and
 the Hilbert space $\mH$ has the following
Hodge decomposition into closed subspaces:
\begin{equation} \label{Hodge}
     \mH = \nul(\Pi_B)\oplus\clos{\ran(\Gamma^*_B)}
            \oplus\clos{\ran(\Gamma)} \ .
\end{equation}
Moreover, $\nul(\Pi_B)=\nul(\Gamma^*_B)\cap\nul(\Gamma)$
and $\clos{\ran(\Pi_B)}=\clos{\ran(\Gamma^*_B)}\oplus
\clos{\ran(\Gamma)}$. When $B_1=B_2=\I$ these decompositions are
orthogonal, and in general the decompositions are topological.

The spectrum $\sigma(\Pi_B)$ is contained in the double sector
$S_{\omega}$, and the operator $\Pi_B$ satisfies resolvent bounds
$$  \|(\I-\tau\Pi_B)^{-1} \|\leq
\frac{C|\tau|}{\dist(\tau,S_\omega)}$$ for all $\tau \in
\cnum\setminus S_\omega$, where $C=C(\|B_1\|,\|B_2\|,\kappa_1,\kappa_2)$.
Such an operator is  of type $S_\omega$ as defined in
\cite{ADMc,AMcNhol}. 

We now introduce further hypotheses which together with (H1-3)
summarize the properties of operators considered in this paper. These form an inhomogeneous version of hypotheses (H4--8) of \cite{AKMC}.

\begin{enumerate}
\item[(H4)] The Hilbert space is
$\mH=L_2(\Omega;\cnum^N)$, where $\Omega \subset \R^n$ and
$n, N \in \N$. 
Here $\Omega$ is a bi-Lipschitz image of $\Omega' \subset
\R^n$, where $\Omega'$ is a smooth domain
which on the complement of a bounded set coincides with either
the empty set $\emptyset$, the half space $\R^n_+$ or $\R^n$.
(By a {\it domain} we mean a connected open set.)
 \item[(H5)] The operators $B_1$
and $B_2$ denote multiplication by matrix valued functions $B_1,
B_2\in L_\infty(\Omega;\mL(\cnum^N))$.
   \item[(H6)] (Localisation)
For every smooth, bounded $\eta:\R^n \longrightarrow \R$ we have that
$\eta \dom(\Gamma) \subset \dom(\Gamma)$, and the commutator $
      M_\eta=[\Gamma,\eta \I]$
is a multiplication operator. There exists $c>0$ so that
$$|M_\eta(x)| \le c
|\nabla\eta(x)|$$ for all such $\eta$ and for all $x\in\rnum^n$. (This implies that the
same hypotheses hold with $\Gamma$ replaced by $\Gamma^*$.)
\item[(H7)] (Cancellation) There exists $c>0$ such that $$\left|
\int_{\Omega} \Gamma u \right| \le c |B|^{1/2} \|u \| \quad
\text{and} \quad \left| \int_{\Omega} \Gamma^* v \right| \le c
|B|^{1/2} \|v \|$$ for every open ball $B$ centred in $\Omega$,
for all $u\in\dom(\Gamma)$ with compact support in $B \cap \Omega$, and
for all $v\in\dom(\Gamma^*)$ with compact support in $B \cap \Omega$.
\item[(H8)] (Coercivity) There exists $\alpha, \beta, c>0$ such that
$$
\| u\|_{H^{\beta} (\Omega;\C^N)} \leq c \| |\Pi|^{\beta} u\|
\quad\text{and}\quad
\| v\|_{H^{\alpha} (\Omega;\C^N)} \leq c \| |\Pi|^{\alpha} v\|
$$ 
for all $u\in \ran(\Gamma^*) \cap \dom (\Pi^2)$
and $v\in \ran(\Gamma) \cap \dom (\Pi^2)$. 
\end{enumerate}
Here $|\Pi|=\sqrt{\Pi^2}$, and
$H^{\beta} (\Omega;\C^N)$ denotes the Sobolev space of order $\beta$ of
$\C^N$-valued functions on $\Omega$.
\begin{rem} \label{const}

In the following theorem and throughout the rest of the paper,
the notation
$a \approx b$ and $b \lesssim c$, for $a,b,c\ge 0$, means that
there exists $C>0$ that depends only on the hypothesis,
so that $ a/C \le b \le Ca$ and $b \le C c$
respectively. 
\end{rem}

\begin{thm} \label{mainthm}
Consider the  operator $\Pi_B=\Gamma+B_1\Gamma^*B_2$ acting in the
Hilbert space $\mH=L_2(\Omega;\cnum^N)$, where $\{\Gamma, B_1,
B_2\}$ satisfies the hypotheses (H1--8). Then $\Pi_B$ satisfies
the quadratic estimate
\begin{equation}
\int_0^\infty\|\Pi_B(\I+t^2{\Pi_B}^2)^{-1}u\|^2\, t\, dt\ \approx\
\|u\|^2 \label{eq:quad}
\end{equation} for all $u \in
\ran(\Pi_B)\subset L_2(\Omega;\cnum^N)$. 
The
comparability constant implicit in the use of ``$\approx$''
depends only on the parameters quantified above
including the bi-Lipschitz constants
implicit in the definition of $\Omega$, and on $\Omega'$. 
\end{thm}

We defer the proof to Section \ref{proof}.

This result implies  that, for every $\omega<\mu<\frac\pi2$, the operator
$\Pi_B$ has a bounded $S^o_\mu$ holomorphic functional calculus in
$\ran(\Pi_B)\subset L_2(\Omega;\cnum^N)$, where $S^o_\mu$ denotes the interior of $S_\mu$.
More to our purposes, it implies the following result. See \cite[Section 2]{AKMC}
for further discussion and proofs.

\begin{cor} \label{cortomainthm22}
Assume the hypotheses of Theorem \ref{mainthm}. Then
$\dom(\Gamma)\cap\dom(\Gamma^*_B)=\dom(\Pi_B) =
\dom(\sqrt{{\Pi_B}^2})$ with
$$\|\Gamma u\|+\|\Gamma^*_Bu\|\approx\|\Pi_Bu\| \approx
\|\sqrt{{\Pi_B}^2}u\|\,.
$$
\end{cor}

\begin{rem} This is equivalent to the statement that 
there is a (non--orthogonal) spectral decomposition
$$\mH=\nul(\Pi_B)\oplus \mH_B^+ \oplus \mH_B^-
$$
 into spectral
subspaces of $\Pi_B$ corresponding to $\{0\}$,
$S_{\omega+}\setminus\{0\}$ and $S_{\omega-}\setminus\{0\}$.
\end{rem}

\subsection{Sobolev spaces}\label{secinter}

We take this opportunity to state some interpolation, trace and extension results
for Sobolev spaces that we need in the next section.

Recall the complex interpolation method.
Let $X \supset Y$ be
Hilbert spaces and $S= \{ z \in \C : 0 < \re z < 1 \}$. Let $\hol
[X,Y]$ denote the Banach space of bounded continuous functions
$f:\overline S \longrightarrow X$ holomorphic on $S$ with $f(z)\in
X$ if $\re z=0$ and $f(z) \in Y$ if $\re z =1$. Each
\emph{interpolation space} $[X,Y]_\theta$, $0<\theta<1$, is given
by
$$[X,Y]_\theta = \{ u(\theta) \in X : u \in \hol [X,Y] \} $$
and inherits a Hilbert space topology from the quotient $\hol
[X,Y]/\{u:u(\theta)=0\}$. A collection of Hilbert spaces
$\{X_s\}_{s \in I}$, where $I \subset \R$ is an
interval, is an \emph{interpolation family} if $X_{(1-\theta)t_1+\theta
t_2} = [X_{t_1},X_{t_2}]_\theta$ for every $0<\theta<1$ and $
t_1,t_2 \in I$.
A good reference for complex interpolation spaces is \cite{lionsb}.

Let $\Omega$, $\Sigma$, $\Sigma_1$ and $\Omega'$ be as described in the
introduction. 

For each $s\in \R$, let
$H^s(\Omega;\C^m)$ denote the fractional order Sobolev space  of
order $s$ of $\C^m$-valued functions on $\Omega$. 
For each $-1 \le s \le 1$, let
$H^s(\Sigma ; \C^m)$ denote the fractional order Sobolev space of
order $s$ of $\C^m$-valued functions on $\Sigma$. 
This can be
defined through localisation arguments by utilizing bi-Lipschitz
parameterizations of $\Sigma$.  
The spaces $H^s(\Sigma ; \C^m)$, $-1 \le s \le 1$, form an interpolation family, 
as do the closed  subspaces
$$
H^s _0(\clos{\Sigma_1} ; \C^m) = \{ u \in H^s(\Sigma ; \C^m) : \supp u \subset
\clos{\Sigma_1} \}\,.$$
More generally, $H^s _0(\clos{\Sigma_1}; \C^m)$ interpolate whenever 
$\Sigma\setminus\clos{\Sigma_1}$ is an extension domain of $\Sigma$. 
By $\Sigma\setminus\clos{\Sigma_1}$ being an extension domain of $\Sigma$,
we mean that the map $R$ restricting distributions in $H^{-1}(\Sigma; \C^m)$ to 
$\Sigma\setminus\clos{\Sigma_1}$ has a right inverse $E$, such that 
$ER: H^s(\Sigma; \C^m)\rightarrow H^s(\Sigma; \C^m)$ is bounded for all $|s|\le 1$.

In the case when $\Omega$ is smooth, i.e. when $\Omega=\Omega'$,
we also make use of the following facts.
\begin{itemize}
\item
  The trace operator $\gamma$ is a bounded map
  $\gamma:H^{s}(\Omega; \C^m) \rightarrow H^{s-1/2} (\Sigma ; \C^m)$
  for $1/2<s<3/2$. 
  There is a bounded extension operator 
  $E: H^{s-1/2} (\Sigma ; \C^m)\rightarrow H^{s}(\Omega; \C^m)$ for
  $0\le s<3/2$ 
  which satisfies $\gamma E=I$  on $H^{s-1/2} (\Sigma ; \C^m)$ when $1/2<s<3/2$.
  In the case $\Omega=\R^n_+$, this map can be contructed as
  in \cite[Section 2.8, Theorem 1b]{aron-smith:bessel}
  with $\alpha=1$.
  \item
  $H^s(\Omega;\C^m)$, $0\leq s < 3/2$, is an interpolation 
  family, as is $
  H^s_0(\Omega ; \C^m)= \sett{u\in H^s(\Omega;\C^m)}{\gamma u=0}
$,  $1/2< s < 3/2$.
  \item
  $H^s_0(\Omega ; \C^m)=[L_2(\Omega ; \C^m),H^1_0(\Omega ; \C^m)]_s$
  when $1/2< s < 1$.
\end{itemize}

%
%  SECTION 3 ##############################################################################
%
%
\section{A Kato square root estimate for systems on domains} \label{slions}

Let us now state a theorem which is somewhat more general than Theorem 
\ref{lions}.
We shall then prove it is a consequence of Corollary \ref{cortomainthm22}, and thus  of Theorem \ref{mainthm}.  Later, in Section \ref{proof}, we shall prove Theorem \ref{mainthm}.

{\bf  Assumptions on $\Omega$, $\mV$ and $S$.} For the remainder of this section, $n, m \in \N$,
 $\Omega$ is an open subset of $\R^n$ which satisfies hypothesis (H4) and has boundary $\Sigma$, and $\mV$ is a closed subspace of $H^1 (\Omega ; \C^m)$
given by \begin{equation} \label{vdef} 
  \mV = \left\{ u \in H^1 (\Omega ; \C^m) :
  \gamma u \in B ^{1/2} ( \Sigma ; \C^m) \right\}\,,
\end{equation}
where
 $B^s( \Sigma ; \C^m)$, $-\frac12 \le s < 1$, is a complex interpolation family of closed subspaces of $H^s(\Sigma ; \C^m)$, and
$B^{1/2} (\Sigma ; \C^m)$ has the following localisation property: 
 whenever $g \in B^{1/2} ( \Sigma;\C^m)$ and 
$\eta: \overline\Omega \longrightarrow \R$ is compactly supported and
Lipschitz, then 
$\eta g \in B^{1/2} (\Sigma ; \C^m)$ with $\|\eta g\|_{B^{1/2}}
\leq c (\|\nabla \eta\|_\infty +\| \eta\|_\infty) \|g\|_{B^{1/2}}$ for some $c$
(independent of $g$ and $\eta$).  In the case when $\Omega=\R^n$, then
$\mV=H^1 (\Omega ; \C^m)$.

Further, $S$ denotes the unbounded operator 
\begin{equation*}
S=\begin{bmatrix}
\I \\
\nabla
\end{bmatrix}
: \dom(S) \subset L_2(\Omega;\cnum^m) \longrightarrow L_2(\Omega;\cnum^{m+nm}) \end{equation*}
with dense domain $\dom(S)=\mV$, and $S^*$ is its adjoint:
\begin{equation*}
S^*=\begin{bmatrix} \I & -
\divv
\end{bmatrix}
: \dom(S^*) \subset L_2(\Omega;\cnum^{m+nm}) \longrightarrow L_2(\Omega;\cnum^m) \,.
\end{equation*}
Then $S$ and $S^*$ are closed and densely defined operators
with $\nul(S)=\{0\}$, $\ran(S)$ closed in $L_2(\Omega;\cnum^{m+nm})$,
 and $\ran(S^*)=L_2(\Omega;\cnum^m)$.

{\bf  Assumptions on $A_1$ and $A_2$.}
Assume that $A_1 \in  L_\infty (\Omega; \mL(\C^m))$ and $A_2
\in L_\infty (\Omega; \mL(\C^{m+nm}))$ satisfy,  for some 
$\kappa_1, \kappa_2>0$,  the accretivity conditions 
 \begin{equation} \label{ACR} \begin{split}   \re (A_1 v, v) &\geq \kappa_1 
\|v\|^2\quad \text{for all} \quad v\in L_2(\Omega;\cnum^m)\, , \\  \re (A_2 S
u, Su) &\geq \kappa_2 \|Su\|^2\quad \text{for all} \quad
u\in\mV\, .
\end{split} \end{equation}
Set    $\omega
:=\tfrac12(\omega_1+\omega_2)$ where \begin{align*}  \omega_1 &:=
\sup\{|\arg(A_1 v, v)|\,:\,v\in L_2(\Omega;\cnum^m), v\neq0\} <\tfrac\pi2
\quad \text{and}\\
     \omega_2 &:= \sup \{|\arg(A_2 S u, S u)| \,:
     \,u\in\mV\setminus\{0\}\}<\tfrac\pi2\,.\end{align*} 

\begin{thm} \label{local}
Suppose that $\Omega$, $\mV$, $S$, $A_1$ and $A_2$ satisfy the above assumptions.
Let $L_A=A_1S^*A_2S$ denote the unbounded operator in 
$L_2(\Omega;\cnum^m)$ with domain $\dom(L_A)=\{u\in\mV:A_2Su\in \dom(S^*)\}$. Then $\sigma(L_A)\subset S_{2\omega+}$ and $L_A$ satisfies resolvent bounds
$\|(\I-\tau L_A)^{-1}\|\lesssim
\frac{|\tau|}{\dist(\tau,S_{2\omega+})}$ for all $\tau \in
\cnum\setminus S_{2\omega+}$,
 so that $L_A$ has a square root $\sqrt{L_A}$ with $\sigma
(\sqrt{L_A})\subset S_{\omega+}$. 

This square root has the Kato square root property 
$\dom(\sqrt{L_A})=\mV$ with 
\begin{equation}
     \|\sqrt{L_A} u\| \approx \|Su\|\approx\| \nabla u \| + \|u \|\end{equation} for all $u\in \mV$. The
comparability constant implicit in the use of ``$\approx$''
depends on $m, c, \|A_1\|_\infty,
\|A_2\|_\infty, \kappa_1, \kappa_2$, and on $\Omega'$,
the bi-Lipschitz constant
implicit in the definition of $\Omega$, and 
constants of interpolation for $B^{s}(\Sigma;\C^m)$.
\end{thm}

We first deduce  that Theorem~\ref{lions} is a consequence of this one.

\begin{proof}[Proof of Theorem~\ref{lions}]
Apply Theorem~\ref{local} with $m=1$, $A_1=I$, $A_2=A$ and 
$B^s(\Sigma;\C)=H^s _0 (\clos\Sigma_1;\C)$, noting that these spaces satisfy the above hypotheses, and that the  sesquilinear form defined in the introduction is
$J_A[u,v]=(ASu,Sv)$, $u,v\in\mV$, with  associated operator $L_A=S^*AS$. 
\end{proof}

We now express Theorem \ref{local}
 in terms of the first order systems presented in 
Section \ref{statement}.
Consider the following operators
$$
     \Gamma=
     \begin{bmatrix}
        0 & 0 \\
        S & 0
     \end{bmatrix},
     \quad
     \Gamma^*=
     \begin{bmatrix}
        0 & S^* \\
        0 & 0
     \end{bmatrix},
     \quad
     B_1=
     \begin{bmatrix}
        A_1 & 0 \\
        0 & 0
     \end{bmatrix},
     \quad
     B_2=
     \begin{bmatrix}
        0 & 0 \\
        0 & A_2
     \end{bmatrix}
$$
in the Hilbert space
$\mH = L_2( \Omega ;\C^m) \oplus  L_2( \Omega ;\C^{m+nm})$. 
They satisfy hypotheses (H1--3), and so have the properties listed in Section 
\ref{statement}.  Moreover 
$$ \Gamma_B^*=B_1\Gamma^*B_2=
     \begin{bmatrix}
        0 & A_1S^*A_2 \\
        0 & 0
     \end{bmatrix}\ ,\ 
     \Pi_B=\Gamma+\Gamma^*_B=
     \begin{bmatrix}
        0 & A_1S^*A_2 \\
        S & 0
     \end{bmatrix} \quad \text{and}$$
$$     {\Pi_B}^2=
     \begin{bmatrix}
        A_1S^*A_2S & 0 \\
        0 & SA_1S^*A_2
     \end{bmatrix}=\begin{bmatrix}
        L_A & 0 \\
        0 & SA_1S^*A_2
     \end{bmatrix}\,.
$$
We remark that $\Gamma,\Gamma^*, \Gamma^*_B, \Pi_B$ and ${\Pi_B}^2$
all have closed range, and that
\begin{gather*}
     \ran(\Gamma) \subset   L_2( \Omega ;
\C^{m+nm}) =
\nul(\Gamma)\quad\text{and}
\\
     \ran(\Gamma^*)= \ran(\Gamma^*_B) =L_2( \Omega ;\C^m)
     \subset \nul(\Gamma^*),\, \nul(\Gamma^*_B).
\end{gather*}
Moreover $\sigma(\Pi_B)\subset S_\omega$ and  $\sigma({\Pi_B}^2)\subset S_{2\omega+}$ with resolvent bounds 
$\|(\I-\tau{\Pi_B})^{-1}\|\lesssim
\frac{|\tau|}{\dist(\tau,S_\omega)}$  and
$\|(\I-\tau^2{\Pi_B}^2)^{-1}\|\lesssim
\frac{|\tau^2|}{\dist(\tau^2,S_{2\omega+})}$ for all $\tau \in
\cnum\setminus S_\omega$.

\begin{prop}  \label{systems}
 Under the above assumptions, the operator $\Pi_B$ satisfies the
quadratic estimate (\ref{eq:quad}) for all $u \in \ran
(\Pi_B)$, and thus $\dom(\Gamma)\cap\dom(\Gamma^*_B)=
\dom(\Pi_B) =
\dom(\sqrt{{\Pi_B}^2})$ with
$$\|\Gamma u\|+\|\Gamma^*_Bu\|\approx\|\Pi_Bu\| \approx
\|\sqrt{{\Pi_B}^2}u\|\,.
$$\end{prop}

We first deduce that Theorem \ref{local} and hence Theorem \ref{lions} is a consequence of Proposition \ref{systems}:

\begin{proof}[Proof of Theorem~\ref{local}]
On restricting the above result to $u\in L_2( \Omega ;\C^m)$ we conclude that
$\dom(\sqrt{L_A})=\dom(S)=\mV$ with the Kato square root
estimate
$$
     \|\sqrt{L_A}u\| \approx \|Su\|=\|\nabla u\|+\|u\|$$ for all $u\in \mV$.
\end{proof}

\begin{rem} It is also a consequence of the quadratic estimate
 (\ref{eq:quad})  that $\Pi_B$ 
 has a bounded $S^o_{\mu}$  
holomorphic functional calculus in $\ran(\Pi_B)\subset\mH$ for 
$\omega<\mu<\pi/2$. Therefore $L_A$ has a  bounded $S^o_{2\mu+}$  
holomorphic functional calculus in $L_2( \Omega ;\C^m)$. This is a generalisation of results in \cite{McN} and \cite{DO}.
\end{rem}

Our task now is to prove  Proposition~\ref{systems}. We do this in two stages. In the first, we show that when the domain $\Omega$ is smooth, then hypotheses (H4--8) are satisfied, and so Theorem \ref{mainthm} and Corollary 
\ref{cortomainthm22} can be applied. This in itself is a new result. In the second stage, we show that the full result is a consequence of the result for smooth domains. 

\begin{proof}[Proof of Proposition~\ref{systems} when $\Omega=\Omega'$]
 Our aim is to verify that $\{\Gamma, B_1, B_2\}$
satisfies hypotheses (H4--8) in $\mH = L_2( \Omega ;\C^m) \oplus  L_2( \Omega ;\C^{m+nm})=
L_2(\Omega;\C^N)$ with $N=2m+nm$.
\begin{itemize}
\item
  Hypothesis (H4) is already assumed, while (H5) follows immediately from the assumptions on $A_1$ and $A_2$.
  \item
  The localisation hypothesis (H6) follows directly from the definition
  (\ref{vdef}) of $\mV$, and the fact that $B^{1/2} (\Sigma ; \C^m)$ 
  satisfies a localisation property.
  \item
  Hypothesis (H7) follows from the fact that 
$\int_{\Omega}\nabla u=0$ and $ \int_{\Omega} \divv v=0$
for $u$ and $v$ with compact support in $\Omega$,
and the use of the Cauchy--Schwarz inequality on the zero order terms.
\item
  To prove (H8), first assume $u\in\ran(\Gamma^*)\cap\dom(\Gamma)$.
Then
$$
  \|u \|_{H^1 (\Omega;\C^N)}= \|\Gamma u\|= \|\Pi u\|=\||\Pi|u\|,
$$
so we can choose $\beta=1$.
Next assume $v=\Gamma u\in\ran(\Gamma)\cap\dom(\Pi^2)$, where $u\in L_2(\Omega;\cnum^m)$.
From Proposition~\ref{regprop} below it follows that
$$
  \|v \|_{H^\alpha (\Omega;\C^N)}=
  \|u \|_{H^{1+\alpha} (\Omega;\C^N)}\lesssim \||\Pi|^{1+\alpha} u\|
  =\| |\Pi|^\alpha \Pi u\|= \| |\Pi|^\alpha v\|,
$$
since $L^{1/2}=|\Pi|$ on $L_2(\Omega;\cnum^m)$.
\end{itemize}
The result now follows on applying Corollary 
\ref{cortomainthm22}.
\end{proof}

We are left with the task of proving the following result.

\begin{prop}    \label{regprop} Under the above assumptions on $\mV$ and $S$, and the smoothness assumption $\Omega=\Omega'$,
  consider the positive operator $L=S^*S$ in $L_2(\Omega;\cnum^m)$.
  Then there exists $\alpha>0$ such that 
$\dom(L^{(1+\alpha)/2})\subset H^{1+\alpha} (\Omega;\C^m)$ with
$$
  \|u\|_{H^{1+\alpha} (\Omega;\C^m)}\lesssim \|L^{(1+\alpha)/2} u \|.
$$
\end{prop}

\begin{rem} Note that $L=-\Delta+I$ with $\dom(L)=\{u\in \mV\,:\, Su\in
\dom(S^*)\}$.
When $m=1$, $\mV$ is defined as in the introduction, and the boundary of $\Sigma_1$ in $\Sigma$ is smooth, then this result can  be derived for any $0<\alpha<\frac12$ as a consequence of results on
mixed boundary value problems
proved by A.~Pryde in \cite{prydethesis}. The proof of  Lemma \ref{inter} is an adaptation of an interpolation argument in \cite{prydethesis}. 
\end{rem}

We prove Proposition~\ref{regprop} with an interpolation and
duality argument using the family
$$ \mV^s = \{ u \in H^s (\Omega; \C^m ) : \gamma u \in B ^{s-1/2} (
 \Sigma ; \C^m) \}$$ for $1/2<s < 3/2$
and the following three lemmas.

\begin{lem} \label{inter}
If $ 1/2 < s \le 1$, then $\mV^{s} \subset \dom(L^{s/2})$ 
with $\|L^{s/2}u\|\lesssim \|u\|_{H^s}$.
\end{lem}

\begin{proof}
By \cite[VI -- Theorem 2.23]{kato}  we have $\mV=\mV^1=\dom(L^{1/2})$. 
It follows that 
$$
  [L_2(\Omega;\cnum^m), \mV]_s = [L_2(\Omega;\cnum^m), \dom(L^{1/2})]_s = \dom(L^{s/2}).
$$

Thus it suffices to prove that $\mV^s \subset [L_2(\Omega;\cnum^m), \mV]_s$.
To this end, let $u\in \mV^s$.
It suffices to show
that there exists $F \in \hol[L_2(\Omega;\cnum^m), \mV]$ with $F (s)=u$. By
definition $g = \gamma u \in B ^{s-1/2} (\Sigma ; \C^m)$.
Therefore there exists
$$G \in \hol [ B ^{-1/2} ( \Sigma ; \C^m), B ^{1/2}
(\Sigma ; \C^m )] \subset \hol[ H ^{-1/2} (
\Sigma ; \C^m), H ^{1/2} (\Sigma ; \C^m)]$$ with $G(s)=g$. 
Let $F_1=EG$ where $E$ is the extension operator mentioned in Section \ref{secinter},
and note that $F_1 \in \hol [ L_2(\Omega;\cnum^m) ,\mV]$. 
Since
$\gamma F_1(s)=\gamma EG(s)=g=\gamma u$, it follows that 
$u-F_1(s) \in H_0 ^s(\Omega; \C^m)$.
Therefore, there exists 
$$
F_2 \in \hol [L_2(\Omega;\cnum^m), H_0 ^1 (\Omega; \C^m)] \subset \hol
[L_2(\Omega;\cnum^m),\mV]
$$
with $F_2(s)=u-F_1(s)$.
Thus
$F=F_1+F_2 \in \hol [L_2(\Omega;\cnum^m), \mV]$ 
and $F(s)=u$. 
This completes the proof.
\end{proof}

\begin{lem}\label{interp} The spaces 
$\{\mV^s\}_{1/2<s<3/2}$ form an interpolation family.
\end{lem}
\begin{proof}
We need to prove that 
$\mV^s=[\mV^{t_1},\mV^{t_2}]_\theta$ for $1/2<t_1<s<t_2<3/2$ and
$s=(1-\theta)t_1+\theta t_2$.
The inclusion $\subset$ is proved in the same way as in 
the proof of Lemma~\ref{inter}.
To prove the incusion $\supset$, let $u=F(\theta)$, where
$$
  F\in \hol[\mV^{t_1},\mV^{t_2}]\subset \hol[H^{t_1}(\Omega;\C^m),H^{t_2}(\Omega;\C^m)].
$$
Thus $u\in H^s(\Omega;\C^m)$.
Furthermore, since 
$\gamma F\in \hol[B^{{t_1}-1/2}(\Sigma;\C^m),B^{{t_2}-1/2}(\Sigma;\C^m)]$
we obtain from the interpolation assumptions on $B^s(\Sigma;\C^m)$ that 
$\gamma u \in B^{s-1/2}(\Sigma;\C^m)$.
This proves that $u\in \mV^s$ and completes the proof. 
\end{proof}

\begin{lem} \label{lemm:sma}
There exists $0<c_0<1/2$ that depends only on $m$ and the constants implicit
in the definition of $\Omega$ such that
for $|\alpha|<c_0$, the form $J:\mV\times \mV\rightarrow \C$
extends to a duality 
$
 J:\mV^{1+\alpha}\times \mV^{1-\alpha}\longrightarrow \C.
$
In particular, we have the estimate 
$$
  \|u\|_{H^{1+\alpha}}\lesssim 
  \sup_{v\in \mV^{1-\alpha}}\frac{|J[u,v]|}{\|v\|_{H^{1-\alpha}}}\ ,
  \quad u\in \mV^{1+\alpha}.
$$
\end{lem}

\begin{proof}
Using the fact that 
$H^\alpha(\Omega;\C^m)$ is the
dual of $H^{-\alpha}(\Omega ; \C^m)$ when $-1/2<\alpha<1/2$
we get
$$
  |J[u,v]|=|(Su,Sv)|\lesssim \|Su\|_{H^\alpha(\Omega ; \C^{m+nm})} 
  \|Sv\|_{H^{-\alpha}(\Omega ; \C^{m+nm})}
  \lesssim\|u\|_{H^{1+\alpha}}\|v\|_{H^{1-\alpha}}\,.
$$
Thus we have an associated bounded operator
$L_\alpha:\mV^{1+\alpha} \rightarrow \left(\mV^{1-\alpha} \right)':
u\mapsto J[u,\cdot]$
when $|\alpha|<1/2$ which is invertible for $\alpha=0$.
By Lemma \ref{interp} and 
the stability result of \v Sne\u\i berg \cite{sneiberg},
there exists a constant $c_0>0$ such that 
$L_\alpha:\mV^{1+\alpha} \longrightarrow \left(\mV^{1-\alpha} \right)'$ 
is an isomorphism when $|\alpha| <c_0$, which proves the lemma.
\end{proof}

\begin{proof}[Proof of Proposition~\ref{regprop}]
Let $c_0$ be the constant from Lemma~\ref{lemm:sma} and let
$0<\alpha<c_0$.
For $u\in\dom( L^{(1+\alpha)/2} )$ we get from 
Lemma~\ref{lemm:sma} and \ref{inter},  that
\begin{align*} 
\| u \|_{H^{1+\alpha}} \lesssim
\sup_{v \in \mV^{1-\alpha}} \frac{|J[u,v]|}{\| v\|_{H^{1-\alpha}}}
&\lesssim \sup_{v \in \dom({L^\frac{1-\alpha}{2}})} \frac{ | (
L^{\frac{1+\alpha}{2}} u , L^{\frac{1-\alpha}{2}} v  )| }{\|
L^{\frac{1-\alpha}{2}} v  \| }\\ &\lesssim\sup_{w \in L_2(\Omega;\cnum^m)} \frac{  |(
L^{\frac{1+\alpha}{2}} u , w )| }{\| w \| } =
\|L^{\frac{1+\alpha}{2}} u \|.
\end{align*}
This completes the proof.
\end{proof}

We have now completed the proof of Proposition \ref{systems} in the case of smooth domains. It remains for us to consider bi-Lipschitz images of smooth domains.  In doing so, we use the following operator theoretic lemma. The proof  is straightforward and we omit it.

\begin{lem}  \label{intertwlemma}
  Let $T:\mH\rightarrow \mH'$ be an isomorphism between Hilbert
spaces, let $\{\Gamma, B_1, B_2\}$ be operators in $\mH$ satisfying
(H1--3).
Assume that $\Gamma'$ satisfies (H1) in $\mH'$
and that $\Gamma' T= T \Gamma$ with 
$\dom(\Gamma') = T \dom(\Gamma)$.

Then $\Pi_B = T^{-1} \Pi'_{B'} T$ with 
$\dom(\Pi'_{B'})= T \dom(\Pi_B)$,
where $\Pi'_{B'}=
\Gamma' +B_1'(\Gamma')^* B_2'$,
$$
  B_1':= TB_1T^*,
  \qquad
  B_2':= (T^{-1})^*B_2T^{-1},
$$
and $\{\Gamma',B_1',B_2'\}$
satisfies (H1--3).
Consequently, if $\Pi'_{B'}$ satisfies quadratic 
estimates, then so does $\Pi_B$.
\end{lem}

\begin{proof}[Proof of Proposition~\ref{systems}]

Suppose that $\Omega, \mV$ and $S$ have the properties specified at the 
beginning of this section, and denote the bi-Lipschitz map from the smooth domain $\Omega'$ with boundary $\Sigma'$ to the
domain $\Omega$ with boundary $\Sigma$ by
 $\rho:\overline{\Omega'}\rightarrow\overline\Omega$.
 
 The map $\rho_0^*$ defined by $\rho_0^* u=u\circ\rho$ is an isomorphism from 
 $L_2(\Omega; \C^m)$ to $L_2(\Omega'; \C^m)$, from
 $H^{1}(\Omega; \C^m)$ to $H^{1}(\Omega'; \C^m)$ and from
$H^{s}(\Sigma; \C^m)$ to $H^{s}(\Sigma'; \C^m)$ 
when $|s|\le 1$, and it commutes with the trace map $\gamma$. On defining
$\mV'=\rho_0^*(\mV)$, we deduce that $\mV$ satisfies the same assumptions on 
 $\Omega'$ as $\mV$ does on $\Omega$.
 Next define $S'$ to be  the unbounded operator 
\begin{equation*}
S'=\begin{bmatrix}
I \\
\nabla
\end{bmatrix}
: \dom(S') \subset L_2(\Omega';\cnum^m) \longrightarrow L_2(\Omega';\cnum^{m+nm}) \end{equation*}
with dense domain $\dom(S')=\mV'$, and let 
$
     \Gamma'=
     \begin{bmatrix}
        0 & 0 \\
        S' & 0
     \end{bmatrix}$.

The operator $\rho_{01}^*=\begin{bmatrix}
        \rho_0^* & 0 \\
        0 & \rho_1^*
     \end{bmatrix}$
     is an isomorphism from $L_2(\Omega;\C^{m+nm})$ to
     $L_2(\Omega';\C^{m+nm})$, where  $\rho_1^*$ denotes the pullback
     $\rho_1^* v:= (d\rho)^t v\circ\rho: L_2(\Omega;\C^{nm})
     \longrightarrow L_2(\Omega';\C^{nm})$. By the chain rule, 
     $S'\rho_0^*=\rho_{01}^*S$.

We can apply the above lemma with
$\mH= L_2(\Omega; \C^m)\oplus L_2(\Omega;\C^{m+nm})$,
$\mH'= L_2(\Omega'; \C^m)\oplus L_2(\Omega';\C^{m+nm})$ and 
$T=\begin{bmatrix}
        \rho_0^* & 0 \\
        0 & \rho_{01}^*
     \end{bmatrix}$, as $\Gamma' T=T\Gamma$.
Now $\Omega',\mV',S',B_1'= TB_1T^*$ and $B_2'= (T^{-1})^*B_2T^{-1}$ satisfy the hypotheses of Proposition \ref{systems}, and we have already proved that $\Pi'_{B'}=
\Gamma' +B_1'(\Gamma')^* B_2'$ satisfies the
quadratic estimate (\ref{eq:quad}) on $\ran(\Pi'_{B'})$. Thus
$\Pi_{B}$ satisfies the
quadratic estimate (\ref{eq:quad}) on $\ran(\Pi_{B})$ as required.
\end{proof}

%
%  SECTION 4 ##############################################################################
%
%
\section{Proof of Theorem \ref{mainthm}} \label{proof}

The proof here is an adaption of
our previous work in \cite{AKMC}.
The main novelty is the inhomogeneity in hypotheses (H7--8).

\begin{defn}
\label{multiop} Define bounded operators in $\mH$  for each $t
\in \R$ by
\begin{align*}
R_t^B&= (\I+it\Pi_B)^{-1}\,,\\P_t^B&=(\I+t^2{\Pi_B}^2)^{-1}=\tfrac
   12(R_t^B+R_{-t}^B)=R_t^BR_{-t}^B\quad \text{and}\\
Q_t^B&=t\Pi_B(\I+t^2{\Pi_B}^2)^{-1}=\tfrac
     1{2i}(-R_t^B+R_{-t}^B)\, , \\
\Theta^B _t&=t\Gamma^*_B(\I+t^2 {\Pi_B}^2 )^{-1}
\end{align*}
In the unperturbed case $B_1=B_2=I$, we write $R_t$, $P_t$ and
$Q_t$ for $R_t ^B$, $P_t ^B$ and $Q_t ^B$.
\end{defn}

To prove Theorem \ref{mainthm} it suffices by the Hodge decomposition
(\ref{Hodge}) and duality considerations as in
 \cite[Proposition 4.8]{AKMC} to
 prove that the square function estimate
\begin{equation} \label{thesqesto}
     \int_0^\infty\|\Theta^B _t P_t u\|^2\,\frac{dt}t\lesssim \|u\|^2
\end{equation}
holds for every $u \in \ran(\Gamma)$ under the hypotheses
(H1--8) stated in Section \ref{statement},
together with the three similar estimates obtained on replacing
$\{\Gamma,B_1, B_2\}$ by $\{\Gamma^* ,B_2, B_1\}$,
$\{\Gamma^*,B_2^*, B_1^*\}$ and $\{\Gamma,B_1^*, B_2^*\}$. As the hypotheses are preserved under these replacements, it suffices to consider 
(\ref{thesqesto}).

We now introduce a \emph{dyadic decomposition} $\triangle$ of
$\Omega$ that is better suited to our circumstance than the
standard dyadic decomposition. It can easily
be constructed using hypothesis (H4).
The decomposition is given by
$\triangle= \bigcup_{j \le j_0} \triangle_{2^j}$ for some $j_0\leq 0$, where 
each $\triangle_{2^j}$ is a collection of Borel subsets $Q$ of $\Omega$
(each of which we refer to as  a \emph{dyadic cube}) such that the
following holds.

\begin{itemize}
\item We have $\Omega=\bigcup_{Q \in \triangle_{2^j}} Q $ for
every integer $j \le j_0$.  
\item We have $Q \cap R=\varnothing$ whenever
$Q, R \in \triangle_{2^j}$ with $Q \neq R.$ 
\item If $R \in \triangle_{2^k}$ and $Q\in
\triangle_{2^j}$  for some $k \le j$,
then either $R \subset Q$ or $ R \cap Q = \varnothing$. 
\item There exists $c\geq 1$ such that for for each $j\leq j_0$ and each $Q \in
 \triangle_{2^j}$, the closure of $Q$ is bi-Lipschitz equivalent to a closed ball of radius $2^j$,
  with bi-Lipschitz constants bounded by $c$.
\end{itemize}

Set $t_0:=2^{j_0}\leq 1$, and for
$0<t\le t_0$, let $\dyadic_t:=\dyadic_{2^j}$ when $ 2^{j-1} < t \le 2^{j}$. 
Note that $|Q|\approx t^{n}$, where $|Q|$ denotes the Lebesgue measure of $Q
\in\triangle_t$.
The {\em dyadic averaging operator} $A_t: \mH
\longrightarrow \mH $ is given by
$$
       A_t u(x) = u_{Q(x,t)}= \barint_{\hspace{-5pt}{Q(x,t)}} u(y)\,  dy = \frac{1}{|{Q(x,t)}|} \int_{Q(x,t)}
       u(y)\,  dy
$$
for every $x \in \Omega$ and $0<t\le t_0$, where  ${Q(x,t)}$ is specified by
$x\in{Q(x,t)} \in \dyadic_t$.

\subsection{Estimates for (\ref{thesqesto})}

To prove the square function estimate (\ref{thesqesto}), we begin
by observing that (H8) implies $\|P_t u \| \le \| |\Pi|^\alpha P_t
u \|$ for every $u \in \ran(\Gamma)$, 
and therefore by spectral theory (because $\Pi$
is self-adjoint) that
$$      \int_{t_0}^\infty\|\Theta^B _t P_t u\|^2\,\frac{dt}t \lesssim
\int_{t_0}^\infty\| (t |\Pi|)^\alpha P_t u\|^2\,\frac{dt}{t^{1+2\alpha}}
\lesssim \|u \|^2 \int_{t_0}^\infty \frac{dt}{t^{1+2\alpha}}  \lesssim
\|u\|^2 $$ where $t_0=2^{j_0}$.
 Thus to prove (\ref{thesqesto}) it suffices to show that
\begin{equation} \label{thesqest}
     \int_0^{t_0}\|\Theta^B _t P_t u\|^2\,\frac{dt}t\lesssim \|u\|^2
\end{equation}
for every $u \in\ran(\Gamma)$.

\begin{defn}
By the {\em principal part} of the operator family $\Theta^B _t$
under consideration, we mean the multiplication operators
$\gamma_t$ defined by
$$
       \gamma_t(x)w= (\Theta^B _t w)(x)
$$
for every $ w\in \cnum^N$. Here we view $w$ on the right-hand side
as the constant function defined on $\Omega$
by $w(x)=w$. It will be proven in Corollary \ref{gammaprops} that
$\gamma_t \in L_2^{\text{loc}}(\Omega; \mL(\cnum^N))$.
\end{defn}

To establish (\ref{thesqest}), we estimate each of the following
three terms separately
\begin{equation}
\begin{split} \label{sqfcn2}
       \int_0^{t_0}&\|\Theta^B _t P_t u\|^2 \frac{dt}t \lesssim
       \int_0^{t_0}\|\Theta^B _t P_t u-\gamma_t A_t P_t u\|^2\frac{dt}t
\\
       &+ \int_0^{t_0}\|\gamma_t A_t (P_t-\I) u\|^2\frac{dt}t
     + \int_0^{t_0}\int_{\Omega} |A_t u(x)|^2 |\gamma_t(x)|^2
\frac{dxdt}t
\end{split}
\end{equation}
when $u\in\ran(\Gamma)$.

We estimate the first two terms in Section~\ref{princsec}, and the
last term in Section~\ref{Carlesonsec}. In the next section we
introduce crucial off--diagonal estimates for various operators
involving $\Pi_B$, and also prove local $L_2$ estimates for
$\gamma_t$.

\subsection{Off--diagonal estimates}

We require off--diagonal estimates for the following class of
operators. Denote $\brac x=1+|x|$, and
$\dist(E,F) =\inf\{|x-y|:x\in E,y\in F\}$ for every
$E,F\subset\Omega$.

\begin{prop} \cite[Section 5.1]{AKMC} \label{pseudoloc}
Let $U_t$ be given by $P^B_t$, $Q^B_t$ or
       $\Theta^B _t$ for every $t>0$ (see Definition \ref{multiop}).
       Then for every $M \in \nnum$ there exists $C_M>0$
       (that depends only on $M$ and the hypotheses (H1--8))
       such that
\begin{equation} \label{odn}
       \|U_t u\|_{L_2(E)} \le C_M \brac{\dist (E,F)/t}^{-M}\|u\|
\end{equation}
whenever $E,F \subset \Omega$ are Borel sets, and $u \in \mH$
satisfies $\supp u\subset F$.
\end{prop}

The proof is be omitted as it is
essentially the same as \cite[Proposition 5.2]{AKMC}. The key
hypothesis used in the proof is (H6). A simple consequence   is that
\begin{equation}  \label{ODest}
       \|U_s u\|_{L_2(Q)}\le \sum_{R\in\dyadic_t}\|U_s(\chi_Ru)\|_{L_2(Q)}
       \lesssim
       \sum_{R\in\dyadic_t}\brac{\dist(R,Q)/s}^{-M}\|u\|_{L_2(R)}
\end{equation}
whenever $0<s\leq t$ and $Q\in\dyadic_t$, where $U_s$ is as
specified in Proposition \ref{pseudoloc}. We also note that the
dyadic cubes satisfy
\begin{equation}  \label{separationest}
       \sup_{Q\in\dyadic_t}\sum_{R\in\dyadic_t}
\brac{\dist (R,Q)/t}^{-(n+1)}\lesssim 1
\end{equation}
and therefore, choosing $M \ge n+1$, we see that $U_t$ extends to
an operator $U_t:  L_\infty (\Omega) \longrightarrow
L_2^{\text{loc}}(\Omega)$.
A consequence of the above results with $U_t=\Theta^B_t$ is:

\begin{cor}  \label{gammaprops}
       The functions $\gamma_t\in L_2^{\text{loc}}(\Omega;
\mL(\cnum^N))$ satisfy the boundedness conditions
$$
       \barint_{\hspace{-5pt}Q} |\gamma_t(y)|^2 \, d y \lesssim  1$$
for all $Q\in\dyadic_t$, $0<t\leq t_0$.  Moreover $\|\gamma_t A_t\|\lesssim 1$
uniformly in $t$.
\end{cor}

\subsection{Principal part approximation}  \label{princsec}
     In this section we prove the principal part approximation
$\Theta^B _t\approx\gamma_t$ in the sense that we estimate the
first two terms on the right-hand side of (\ref{sqfcn2}). The
following lemma is used in estimating the first term.

\begin{lem}  \label{poincarelem}
If $0<t\leq t_0$, $Q\in\dyadic_t$ and $M > 2n$, then
we have
$$
      \int_{\Omega}|u(x)-u_Q|^2\brac{\dist(x,Q)/t}^{-M} \,  dx \lesssim
      \int_{\Omega}(|t\nabla u(x)|^2 + |tu(x)|^2)
       \brac{\dist(x,Q)/t}^{2n-M}\,  dx
$$
for every $u$ in the Sobolev space $H^{1} (\Omega; \cnum^N)$.
\end{lem}

\begin{proof}  In the case when $\Omega$ is a smooth domain  one can use reflection techniques to 
construct an extension operator $\mE:H^{1} (\Omega; \cnum^N)\longrightarrow H^{1} (\R^n; \cnum^N)$ such that
$$\int_{\R^n}|\nabla (\mE u)(x)|^2 \brac{\dist(x,Q)/t}^{2n-M}\,  dx
\lesssim \int_{\Omega}(|\nabla u(x)|^2 + |u(x)|^2)
 \brac{\dist(x,Q)/t}^{2n-M}\, dx\,.$$
The desired estimate then follows from the corresponding result on $\R^n$ 
\cite[Lemma 5.4]{AKMC}, noting that the set $Q$ used there does not need to be a Euclidean cube, but merely satisfy $|Q|\approx t^n$.

In the general case of a domain which is bi-Lipschitz equivalent to a smooth domain,  the bi-Lipschitz  parametrization of $\Omega$  gives the required inequality, except for the fact that $u_Q$ is replaced by a constant $c =c(u,Q)$. But this suffices, because $\int_{\Omega}|u_Q-c|^2\brac{\dist(x,Q)/t}^{-M}   dx \lesssim t^n|u_Q-c|^2 \lesssim 
\int_{\Omega}|u(x)-c|^2\brac{\dist(x,Q)/t}^{-M}  dx$.
\end{proof}

We now estimate the first term in the right-hand side of
(\ref{sqfcn2}).

\begin{prop}  \label{secondtermprop}
       For all $u\in \ran(\Gamma)$, we have
$$
       \int_0^{t_0}\|\Theta^B _t P_t u-\gamma_t A_t P_t u\|^2\,\frac{dt}t
       \lesssim  \|u\|^2\, .
$$
\end{prop}
\begin{proof}
Using Proposition~\ref{pseudoloc}, Lemma~\ref{poincarelem} and
estimate (\ref{separationest}) we get, as in \cite{AKMC}, that
\begin{equation*} 
       \|\Theta^B _t v -\gamma_t A_t v  \| \lesssim 
    t \|v\|_{H^{1} (\Omega; \cnum^N)}
\end{equation*}
for every $v$ in the Sobolev space $H^{1} (\Omega; \cnum^N)$. 
Since $\Theta^B _t -\gamma_t A_t$ is bounded on
$\mH$, we have by interpolation and then (H8) that
$$ \|\Theta^B _t v-\gamma_t A_t v  \| \lesssim  t^\alpha \|v\|_{H^\alpha(\Omega;\C^N)}
\lesssim \| (t|\Pi|)^\alpha v \|
$$
for every $v \in \ran(\Gamma) \cap \dom(\Pi^2)$. 

Taking $v=P_tu$, we then have
$$
       \int_0^{t_0}\|\Theta^B _t P_t u-\gamma_t A_t P_t u\|^2\,\frac{dt}t
       \lesssim \int_0^{t_0}\| (t|\Pi|)^\alpha P_t u\|^2\,\frac{dt}t
       \lesssim  \|u\|^2\, .\quad  \qedend
$$
The last inequality above follows from spectral theory. This
completes the proof.
 \end{proof}

We use the following lemma to estimate  the second term on the
right-hand side of (\ref{sqfcn2}), and also in the proof of Lemma
\ref{lem:ncar1} (c.f. Lemma 5.15 of \cite{AHLMcT}).

\begin{lem} \label{interpolationlemma} 
Let $\Upsilon$ be either $\Pi$, $\Gamma$ or $\Gamma^*$. Then we
have the estimate
\begin{equation}  \label{interpolest}
\left| \barint_{\hspace{-5pt} Q} \Upsilon u  \right|^2 \lesssim
\frac 1t \left( \barint_{\hspace{-5pt}Q} |u|^2 \right)^{1/2}
\left( \barint_{\hspace{-5pt}Q} | \Upsilon u|^2 \right)^{1/2} +
 \barint_{\hspace{-5pt}Q} |u|^2
\end{equation}
for all $Q \in\dyadic_t$ and $u \in \dom(\Upsilon)$.
\end{lem}
\begin{proof}
Let $\tau=(\int_Q|u|^2)^{1/2}(\int_Q|\Upsilon u|^2)^{-1/2}$. If $\tau\ge
t$, then (\ref{interpolest}) follows directly from the
Cauchy--Schwarz inequality. If $\tau\leq t$, let $\eta\in
C^\infty_0(Q)$ be a real-valued bump function with
$|\nabla\eta|\lesssim 1 /\tau$ such that $\eta(x)=1$ whenever $x \in
Q$ satisfies $ d(x, \rnum^n\setminus~Q)
\ge  \tau.$ Then by  hypothesis (H7),
the Cauchy--Schwarz inequality, and the fact that
$ \left| \{ x \in Q  :  \dist(x , \R^n \setminus Q) \le \tau
\} \right| \lesssim \tau t^{n-1}$, 
we obtain
\begin{equation*} \begin{split}
       \left|\int_Q \Upsilon u\right| &=\left|
\int_Q (1-\eta)\Upsilon u + \int_Q \eta \Upsilon u\right|
       \\ & =\left| \int_Q (1-\eta)\Upsilon u
       +\int_Q [\eta,\Upsilon ]u+ \int_Q \Upsilon (\eta u) \right| \\
       &\lesssim   (\tau t^{n-1})^{1/2}
\left(\int_Q|\Upsilon u|^2\right)^{1/2} \\ & \quad  +
 \|\nabla\eta\|_{\infty} ( \tau t^{n-1})^{1/2}
\left(\int_Q|u|^2\right)^{1/2} +
 |Q|^{1/2} \left(
\int_{Q} |u|^2 \right)^{1/2}
\end{split} \end{equation*}
which leads to (\ref{interpolest}) on substituting the chosen
value of $\tau$. \qedend \end{proof}

We now estimate the second  term in the right-hand side of
(\ref{sqfcn2}).

\begin{prop}  \label{thirdtermprop}
       For all $u\in \mH$, we have
$$
       \int_0^{t_0} \|\gamma_t A_t (P_t-\I) u\|^2\,\frac{dt}t \lesssim  \|u\|^2\, .
$$
\end{prop}

\begin{proof}
Corollary~\ref{gammaprops} shows that $\|\gamma_t A_t\|\lesssim 1$
and since $A_t^2=A_t$ it suffices to prove the square function
estimate with integrand $\|A_t (P_t-\I) u\|^2$. If $u\in\nul(\Pi)$
then this is zero. If $u\in\ran(\Pi)$ then by spectral
theory we can write $u=2\int_0^\infty Q_s^2u\frac{ds}s$. The
result will follow from a Schur estimate and the spectral theory
estimate $\int_0 ^\infty \|Q_t u \|^2 \, \frac{dt}{t} \le \|u
\|^2$ once we have obtained the bound
$$
       \|A_t (P_t-\I)Q_s\|\lesssim \min\{\tfrac st,\tfrac ts\}^{1/2}$$
for all $s>0$ and $0 <t \le t_0$.

Note that $  (\I-P_t)Q_s = \tfrac ts Q_t(\I-P_s)$
 and $ P_tQ_s =\tfrac st Q_tP_s$
for every $s,t>0$. Thus, if $t \le\min(s,t_0)$, then $$\|A_t
(P_t-\I)Q_s\|\lesssim \|(P_t-\I)Q_s\|\lesssim  t/s\,,$$ while if
$s <t \le t_0$, then $$\|A_t (P_t-\I)Q_s\|\lesssim
\|P_tQ_s\|+\|A_tQ_s\|\lesssim s/t+\|A_tQ_s\|\,.$$ To estimate
$\|A_tQ_s\|$, we use Lemma~\ref{interpolationlemma} with
(\ref{ODest}) and (\ref{separationest}) to obtain
\begin{equation*} \begin{split}
       \|A_tQ_s u\|^2 &= \sum_{Q\in\dyadic_t} |Q|\,
       \bigg| \barint_{\hspace{-5pt}Q} s\Pi(\I+s^2 \Pi ^2 )^{-1}u \bigg|^2 \\
       & \lesssim \frac st\sum_{Q\in\dyadic_t}
       \bigg( \int_{Q}|P_su|^2 \bigg)^{1/2}
       \bigg( \int_{Q}|Q_su|^2 \bigg)^{1/2}  + s^2  \int_Q
       |P_s u |^2 \\
       & \lesssim \frac st\sum_{Q\in\dyadic_t}
         \bigg(
\sum_{R\in\dyadic_t}\brac{d(R,Q)/t}^{-(n+1)}\|u\|_{L_2(R)}
\bigg)^2 + \left(\frac st\right)^2\|u\|^2
     \\
       & \lesssim \frac st \sum_{Q\in\dyadic_t}
         \bigg( \sum_{R'\in\dyadic_t}\brac{d(R',Q)/t}^{-(n+1)} \bigg)
         \bigg( \sum_{R\in\dyadic_t}\brac{d(R,Q)/t}^{-(n+1)}\|u\|_{L_2(R)}^2
         \bigg) \\&+ \left(\frac st\right)^2\|u\|^2 \lesssim \frac st \|u\|^2
\end{split} \end{equation*}
which completes the proof. \qedend \end{proof}

We have now estimated the first two terms in the right-hand side
of (\ref{sqfcn2}).

\subsection{Carleson measure estimate}  \label{Carlesonsec}

In this subsection we estimate the third term in the right-hand
side of (\ref{sqfcn2}). To do this we reduce the problem to a
Carleson measure estimate. Recall that a measure $\mu$ on $\Omega
\times (0,t_0)$ is said to be \emph{Carleson}  if
$\|\mu\|_{\mC}=\sup_{Q\in\triangle}|Q|^{-1}\mu(R_Q)<\infty$.
 Here
and below, $R_Q = Q \times (0,2^j)$ denotes the {\it Carleson box} over $Q \in \triangle_{2^j}$. For such $Q$ we define
$ \lambda Q = \{ x \in \R^n : \dist (x,Q) \le (\lambda-1) 2^j
\}$ when $\lambda \ge 1$.

 We now recall the
following theorem of Carleson.

\begin{thm}\cite[p.\ 59]{stein:harm}
If $\mu$ is a Carleson measure on $\Omega \times (0,t_0)$ then
$$ \iint_{\Omega \times (0,t_0)} |A_t u (x)|^2 \, d \mu(x,t) \le C
\|\mu \|_{\mC} \| u \|^2$$ for every $u \in \mH$. Here $C>0$ is a
constant that depends only on $n$.
\end{thm}

Thus, in order to prove (\ref{sqfcn2}) it suffices to show that
\begin{equation} \label{ncar1} \iint_{R_Q} | \gamma_t (x) |^2 \,
\frac{dxdt}{t} \lesssim |Q|\end{equation} for every dyadic cube $Q
\in \triangle$. 

Define a measure $\nu $ on $\Omega \times (0,t_0)$ 
by $d\nu = \chi (x,t) \frac{dxdt}{t}$, where $\chi$
is the characteristic function defined by $\chi (x,t)=1$ if there
is $x \in Q \in \dyadic_t$ with $4Q \setminus \Omega \neq
\emptyset$; otherwise let $\chi (x,t)=0$. It follows from (H4)
that $\nu$ is a Carleson measure. From Corollary \ref{gammaprops}
we then see that $|\gamma_t(x)|^2 \, d \nu (x,t)$ is a Carleson
measure. Now, the sum of two Carleson measures is again Carleson.
Therefore, to prove (\ref{ncar1}) it remains to consider the case
when $Q \in \triangle$ with $4Q \subset \Omega$.

Fix such a cube $Q$ and set $\sigma>0$; its value to be
chosen later. Let $\mF$ be a finite set consisting of $\nu \in \mL
(\C^N)$ with $|\nu|=1$, such that $\bigcup _{ \nu \in \mF} K_\nu =
\mL ( \C^N) \setminus \{0\},$ where
$$K_\nu = \left\{ \nu' \in \mL(\C^N) \setminus \{0\} : | \tfrac{\nu'}{|\nu'|} - \nu |
\le \sigma  \right\}.$$ To prove (\ref{ncar1}) it suffices to show
that
\begin{equation} \label{ncar2}
\iint_{\substack{ (x,t) \in R_Q \\ \gamma_t(x) \in K_\nu }} |
\gamma_t(x) |^2 \, \frac{dxdt}{t} \lesssim |Q|\end{equation} for
every $\nu \in \mF$. By a standard stopping time argument as used
in \cite[Section 5]{AHLMcT}, in order to prove (\ref{ncar2}) it
suffices to prove the following claim.

\begin{prop} \label{prop:ncar2}
There exists $\beta
>0$ such that for every dyadic
cube $Q\in \triangle$ with $4Q \subset \Omega $,
and for every $\nu \in \mL ( \C^N)$ with $|\nu |=1$, there is a
collection $\{Q_k\}_{k} \subset \triangle$ of disjoint subcubes of
$Q$ such that $ |E_{Q,\nu}|
>  \beta |Q|$ where $E_{Q,\nu }= Q \setminus \bigcup_k Q_k$, and such
that
$$ \iint_{\substack { (x,t) \in  E^* _{Q,\nu} \\ \gamma_t(x) \in
K_\nu }} | \gamma_t(x) |^2 \,\frac{dx dt}{t} \lesssim |Q|$$ where
$E^* _{Q ,\nu } = R_Q \setminus \bigcup_k R_{Q_k}$.
\end{prop}

Let $Q\in\triangle_\tau$ and $\nu$ be as in the above proposition. Choose $\hat w,
w \in \C^N$ with $|\hat w|=|w|=1$ and $\nu^* (\hat w )=w$. Let
$\eta_Q$ be a smooth cut-off function with range $[0,1]$, equal to
$1$ on $2Q$, with support in $4Q$, and such that $\| \nabla \eta_Q
\|_\infty \le 1/\tau$. Define $w_Q =
\eta_Q w $, and for each $\epsilon>0$, let
$$ f_{_{Q,\epsilon}}^w =w_{_Q}-\epsilon \tau
i\Gamma(1+\epsilon \tau i\Pi_B)^{-1}w_{_Q}\\
=\left(1+\epsilon \tau i\Gamma_B^*\right)(1+\epsilon \tau
i\Pi_B)^{-1}w_{_Q}\,.$$

\begin{lem} \label{lem:ncar1}
We have $\|f^w _{Q,\epsilon } \| \lesssim |Q|^{1/2},$
\begin{align*}
\iint_{R_Q} |\Theta_t ^B f^w_{Q,\epsilon}|^2 \, \frac{dxdt}{t}
\lesssim \frac{1}{\epsilon^2} |Q| \quad \text{and}\quad
 \left| \barint_{\hspace{-5pt}Q} f^w_{Q,\epsilon} - w \right| \le C\epsilon^{1/2}\end{align*} for every $\epsilon >0$. Here $C$
is a constant that depends only on hypotheses (H1--8).
\end{lem}

\begin{proof}
The first and second estimates follow as in \cite{AKMC}. To obtain
the last estimate, we use the fact that $\tau \le 1$ and also
Lemma~\ref{interpolationlemma} with $\Upsilon=\Gamma$ and
$u=(\I+\epsilon \tau i\Pi_B)^{-1}w_Q$ to show that
\begin{equation*}\begin{split}
       \left| \barint_{\hspace{-5pt}Q} f^w_{Q,\epsilon} - w\right|
      & = \left|\barint_{\hspace{-5pt}Q} \epsilon \tau\Gamma (\I+\epsilon \tau i\Pi_B)^{-1}w_Q\right|
     \\
       &\lesssim \epsilon^{1/2} \left(\barint_{\hspace{-5pt}Q} |(\I+\epsilon
\tau i\Pi_B)^{-1}w_Q|^2\right)^{1/4}
       \left(\barint_{\hspace{-5pt}Q} | \epsilon \tau \Gamma (\I+\epsilon
\tau i\Pi_B)^{-1}w_Q|^2\right)^{1/4}
\\ & \quad + \epsilon \tau \left( \barint_{\hspace{-5pt}Q}
 |(\I+\epsilon \tau i\Pi_B)^{-1}w_Q|^2\right)^{1/2}
       \lesssim  \epsilon^{1/2}.
\end{split} \end{equation*}
This completes the proof. \qedend \end{proof}

The proof of Proposition~\ref{prop:ncar2} can
now be completed exactly as in \cite{AKMC}.
Therefore the last term in (\ref{sqfcn2}) is bounded by a constant
times $\|u\|^2$.
This proves the square function estimate
(\ref{thesqesto}) and thus Theorem \ref{mainthm}.

\bibliographystyle{acm}

\end{document}